%% file: math148paper.tex
\documentclass[12pt]{amsart}

\usepackage[export]{adjustbox}

\usepackage{amssymb}
\usepackage{bm}
\usepackage[usenames,dvipsnames]{xcolor}
\usepackage{enumerate}
\usepackage{graphicx}
\usepackage{psfrag}
\usepackage{nicefrac}
\usepackage{color}
\usepackage{url}
\usepackage{enumerate}
\usepackage{subcaption}
\usepackage{mathtools}

\usepackage{xy}
\input xy
\xyoption{all}

\newcommand{\excise}[1]{}

\theoremstyle{definition}

\numberwithin{equation}{section}

\newcommand{\ring}[1]{\ensuremath{\mathbb{#1}}}

\renewcommand\>{\rangle}

\newcommand\RR{\ring{R}}

\newcommand\ZZ{\ring{Z}}

\newcommand{\FF}{\mathbb{F}}

\begin{document}

\mbox{}
\title[Discovery learning in an interdisciplinary course on finite fields]{Discovery learning in an interdisciplinary course \\ on finite fields and applications}

\author[C.~O'Neill]{Christopher O'Neill}
\address{Mathematics and Statistics Department\\San Diego State University\\San Diego, CA 92182}
\email{cdoneill@sdsu.edu}

\author[L.~Silverstein]{Lily Silverstein}
\address{Mathematics and Statistics Department\\California State Polytechnic University\\Pomona, CA 91768}
\email{lsilverstein@cpp.edu}

\date{\today}

\begin{abstract}
The authors describe their approach to teaching a course on finite fields and combinatorial applications, including block designs and error-correcting codes, using a hybrid of lectures and active learning.  Under the discussed classroom model, there are two lecture days and two discovery-based discussion days each week.  Discussions center around group activities that build intuition for abstract concepts while avoiding excessive technical machinery.  Teaching this course presents some unique challenges, as much of the content typically appears in a second course in abstract algebra, yet the students exhibit a wide range of mathematical preparation. 
\end{abstract}

\maketitle


\input{sec-introduction}

\input{sec-coursecontent}

\input{sec-coursestructure}

\input{sec-sample-discussions}

\input{sec-2ndtime}
\input{sec-conclusion}

\bibliographystyle{plain}
\bibliography{math148paper.bib}


%



\end{document}

%% file: sec-introduction.tex

\section{Introduction}
\label{sec:intro}

The present work concerns an upper division course on discrete mathematics, with particular attention to block designs, error-correcting codes, and several prerequisite topics in abstract algebra (see Section~\ref{sec:content} for a more thorough content overview).  The~class is usually split roughly evenly between mathematics majors and computer science majors. A~prerequisite is proof-writing experience from a prior class, but the students have highly varied mathematical backgrounds.  Based on informal surveys at the start of each course, about a quarter of the students have already taken abstract algebra, while another quarter have never seen modular arithmetic. This variance means the abstract algebra content must be rigorous enough for junior and senior math majors, while simultaneously communicating enough intuition, motivating examples, and combinatorial applications to hold the interest of computer science majors.

The authors taught this course together in two consecutive years, using a classroom model in which each class meeting is designated as either a ``lecture day'' or a ``discussion day.''  Discussion days center around working in small groups, usually at a chalk- or whiteboard, on discovery learning worksheets with open ended problems designed to facilitate student discovery of new course material.  During lecture days, the instructor sets the pace while introducing new material and preparing students for discussions.  

There is a plethora of documented evidence that active learning benefits students (e.g., \cite{activelearning}), and it is not hard to see why.  
Students take ownership of the material, discovering much of the content on their own rather than simply taking dictation.  There is more opportunity for instructors to observe and interact with the students, and to identify and correct misconceptions as they appear. Open-ended and inquiry-oriented activities are an opportunity to emphasize \textit{process skills} (see, e.g., \cite{hanson-POGIL,POGIL-chem}), in addition to particular content. One of the process skills we emphasized was \textit{mathematical investigation} through tasks like creating conjectures based on examples.

On the other hand, designing and implementing a course based on student inquiry poses new challenges for instructors. In a typical flipped classroom, (e.g., employing the Moore or modified Moore method \cite{mahavier-modified-texas}), the increased workload for both students and teachers can be significant.  Students are responsible for learning most or all of the new material outside of class, requiring a greater time investment than they may be accustomed to, while teachers must create the additional materials, such as instructional videos, that the students are using outside of class. The idea of having to create a semester's worth of instructional videos before attempting a flipped classroom is enough to turn otherwise interested parties away; our blended course design is one way to gain experience with activity-based classrooms with a lower barrier to entry.

Despite the increased planning, care, and time commitment of the instructor, this effort (and its benefits) may be invisible to students. A recent study that compared passive lecture to active learning, in classrooms with identical handouts and assessments, found that students in active classrooms learned more but \textit{perceived} that they learned less \cite{Deslauriers19251}. As junior instructors, we were concerned not only with student perceptions of our teaching/their learning, but also with the perceptions of senior faculty and administration. It can be daunting to implement any unusual or experimental methods in the classroom, out of a fear that the outcomes and evaluations will be subject to extra scrutiny and criticism. Our approach could be a good option for instructors with similar trepidation.

Another concern is maintaining academic standards and sufficiently introducing all course topics, particularly if the course is a prerequisite for other courses. 
By incorporating lectures and carefully choosing which material to deliver via lecture, we found we retained control over the pace of the course. We also retained some flexibility to adjust the ratio of passive lecture to active learning, in case we encountered unforeseen difficulties.  Although a change wasn't necessary in our situation, having this flexibility gave us greater confidence to try unfamiliar methods.

 Additionally, blending lecture with active learning allowed us to introduce the specific definitions and notation that conformed with the other courses at our university. In completely inquiry-based curricula, such as the excellent Teaching Abstract Algebra For Understanding (TAAFU), students are given freedom to develop their own definitions, notation, theorem statements, etc., as they discover algebraic concepts \cite{larsen-design-taafu,larsen-groups-taafu,larsen-instructor-support}. In a university where most courses are based on textbooks and taught via passive lecture, it would be difficult to implement this style for one particular class, then expect those students to rejoin the textbook curriculum without disadvantage.

The rest of this paper goes into detail about our specific implementation.  All course materials are freely available at the following webpages, and some samples are included as figures in Sections~\ref{subsec:sample-discussion-modular} and~\ref{subsec:sample-discussion-projective}.  

\begin{center}
\url{https://cdoneill.sdsu.edu/teaching/w17-148/}

\url{https://cdoneill.sdsu.edu/teaching/w18-148/}
\end{center}


Teaching the same course twice gave us the opportunity to make small adjustments to the course design and content. In Section \ref{sec:2ndtime}, ``Changes the second time," we highlight the most useful of these changes.

In addition to the discussion participation points (see Sections~\ref{subsec:prelim-problems} and~\ref{subsec:enforcing-attendance}), grades are based on weekly homework, one midterm exam, and a final exam.  This conforms to the usual expectations of the university and the students, which, again, is one of the deliberate features of our classroom model.  Like lecture days, exams help set a pace to ensure conformity with the department syllabus.  

\subsection*{Acknowledgements}

The authors would like to thank Stephan Ramon Garcia for his advice and support, and several anonymous referees for their thorough feedback.

%% file: sec-coursecontent.tex

\section{Course content}
\label{sec:content}

\subsection{Content survey}

An \emph{error-correcting code} is a collection of \emph{codewords} used to represent messages in such a way that the recipient can reliably detect and/or correct a minimum number of errors occurring in transmission.  For example, if the message is a series of bits, the code
$$0 \mapsto 000, \qquad 1 \mapsto 111$$
allows the recipient to correct an error where a single bit is transmitted incorrectly.  More specifically, if a transmission is received as
\begin{center}
010 111 110 000,
\end{center}
then the recipient detects incorrect bits in the first and third letters, and correcting each to the nearest valid codeword yields the intended message $0110$.  

The construction of space-efficient error correcting codes is one of the primary goals in coding theory, and has countless applications in computer science and engineering.  One method of constructing highly space-efficient error-correcting codes (including some so-called \emph{perfect} codes) is using \emph{block designs}, which are collections of sets of equal size, called \emph{blocks}, with certain intersection properties.  For example, the collection
\begin{center}
$\begin{array}{l@{\qquad}l@{\qquad}l@{\qquad}l}
\begin{array}{rrr}
\{1, & 2, & 3 \} \\
\{4, & 5, & 6 \} \\
\{7, & 8, & 9 \} \\
\end{array}
&
\begin{array}{rrr}
\{1, & 4, & 7 \} \\
\{2, & 5, & 8 \} \\
\{3, & 6, & 9 \} \\
\end{array}
&
\begin{array}{rrr}
\{1, & 5, & 9 \} \\
\{2, & 6, & 7 \} \\
\{3, & 4, & 8 \} \\
\end{array}
&
\begin{array}{rrr}
\{1, & 6, & 8 \} \\
\{2, & 4, & 9 \} \\
\{3, & 5, & 7 \} \\
\end{array}
\end{array}$
\end{center}
forms a block design in which every pair of elements appears together in exactly one block.  In addition to their use in constructing error-correcting codes, block designs are used in statistics for experiment design and group testing.

Block designs and error-correcting codes are studied in the second half of the course, with a focus on constructions using the algebraic and geometric properties of finite fields.  The first half of the course is spent introducing a collection of carefully chosen prerequisite topics, starting with modular arithmetic, moving on to polynomial rings and factorization, and concluding with finite fields. These topics are developed with enough depth to understand the later material, but with less rigor than in the university's abstract algebra sequence.  Students are introduced to the idea of block designs and error-correcting codes on the first day of class, and throughout the first half of the course they are reminded of these applications as motivation for the more abstract prerequisite topics.  

Some homework problems are inspired by exercises and examples in the optional course textbook \emph{Discrete Mathematics} by Biggs \cite{biggsdiscrete}.  

\subsection{Minimizing technical machinery}

Given the volume of algebraic material to cover and the students' varied algebra backgrounds, difficult vocabulary and technical proofs are carefully avoided when they are unnecessary or unenlightening.  Some concepts, such as \emph{isomorphisms}, do not play a central role in later topics and can simply be conveyed with a couple of well-chosen examples rather than a formal definition.  Other concepts require more care to avoid, such as \emph{Galois groups}, which, though not essential to this course, are used in most proofs of the fundamental theorem of finite fields.  As such, the instructor must carefully choose which statements to prove formally and which to convey using illustrative examples.  

One of the more difficult topics to address is \emph{quotient rings}, which students must use explicitly in order to work with certain finite fields like $\FF_4 = \ZZ_2[z]/\<z^2 + z + 1\>$, the field with 4 elements.  Rather than introduce quotient rings in full generality, they are introduced in a computational manner. For example,
$$\FF_4 = \ZZ_2[z]/\<z^2 + z + 1\> = \{0, 1, z, z + 1\},$$
where the product of two elements is defined as its remainder from polynomial long division by $z^2 + z + 1$, e.g.,
$$(z + 1) \cdot (z + 1) = z^2 + 1 = 1 \cdot (z^2 + z + 1) + z \equiv z \bmod (z^2 + z + 1).$$
Students can verify the field axioms using the fact that $z^2 + z + 1$ is irreducible in $\ZZ_2[z]$. 
All quotient rings encountered in the course are presented this way, i.e., as a~polynomial ring in one variable modulo a single polynomial, so this computational lens is both sufficient for subsequent course material and appealingly concrete for practical use.

%% file: sec-coursestructure.tex

\section{Course structure}
\label{sec:structure}


Our classroom model makes a rigid distinction between lecture days and discussion days.  The class met four days a week, and we usually held lectures on Mondays and Wednesdays, and discussions on Thursdays and Fridays.  We did occasionally deviate, but students were always notified ahead of time so that they knew what style to expect.  
Dedicated discussion days provide a setting for inquiry-based learning, while
 lecture days give instructors more control over the pace of the material than in a fully inquiry-based course.  Discussion problems are hand-picked to be reasonably accessible via discovery, leaving more technical or lengthy topics to be covered in lecture.


One of the primary challenges in planning the course is deciding which content should be introduced during lectures and which should be introduced during discussions.  For some topics, definitions and basic examples are introduced in lecture, leaving deeper connections to be discovered in subsequent discussions.  As an example, one lecture states that the collection of lines in the 2-dimensional vector space $\FF_q^2$ satisfies the definition of a certain type of block design.  Students then discover in the following discussion that for $\FF_q^2$, each block design axiom follows from a geometric property (e.g., two points uniquely determine a line).  

Other topics are first encountered in discussion, where students discover a theorem in practice before seeing it formally.  For example, after working with the field axioms for several weeks, students are asked in discussion to use the axioms to fill in the addition and multiplication tables for all possible fields of small (finite) size.  One of their primary tools is the ``sudoku rule''~\cite{cooksudoku}, which states that every element appears exactly once in each row and column of the addition table, and each nonzero element appears exactly once in each nonzero row and column of the multiplication table (these are consequences of the respective invertibility axioms).  Students soon discover that there is only one way to fill in the $3 \times 3$ and $4 \times 4$ operation tables while satisfying all field axioms, and that it is impossible to do this for the $6 \times 6$ operation tables.  These facts are a consequence of the fundamental theorem of finite fields, which is presented in the next lecture and states that (i)~any finite field of a given size is unique up to a relabeling of its elements, and (ii)~the size of any finite field equals a prime power.

\input{sec-discussions}

%% file: sec-discussions.tex

\subsection{Discovery learning days}
\label{sec:discussions}

Problems completed during the discussion days are often open ended, prompting students to formulate conjectures from examples rather than verifying existing claims.  Students are not required to submit their answers to the discussion questions, freeing them to focus on exploring and discovering the material rather than on getting a good grade on their write-up.  In addition to acquiring experience with the process of mathematical investigation, discovery learning questions give students a feeling of ownership over the material, helping to build confidence in their mathematics skills.  Indeed, as the course progresses, many students become more willing to try out examples and see what happens, without a concern for finding the right answer on the first try.  

Some college students have learned to dread group assignments, due in part to the sometimes unequal distribution of work.  
Another potential issue is unequal distribution of mathematical background:\ better-prepared students are impatient to move through questions quickly, while less-prepared students are frustrated when the conversation moves too fast.  We attempted to avoid these sources of reluctance and frustration by clearly separating in-class work from graded work.  Discussion problems are never turned in for credit, eliminating the pressure to rush through problems or to forgo collaboration in a divide-and-conquer scheme.  Moreover, students are not expected to complete every discussion problem; the instructor can decide ahead of time which parts are the most crucial and guide the groups appropriately to ensure they get to these, leaving the remaining parts as extra practice for groups that move quickly.  


\subsection{Preliminary problems}
\label{subsec:prelim-problems}

The day before each discussion, students are expected to complete a short (at most 10 minutes) preliminary assignment that is computational in nature (i.e., no proofs).  For example, these were the preliminary problems assigned before our first discussion of the quarter.  
\begin{enumerate}[(P1)]
\item 
Fill in the multiplication table for $\mathbb Z_5$. [A table with row and column labels is provided on the worksheet.]

\item 
Find all $x\in\mathbb Z_7$ that satisfy $x^2 = [4]_7$.

\end{enumerate}
The preliminary problems force students to review their notes from the most recent lecture, so they can hit the ground running on discussion days, and also help to tie together the two classroom formats.  Students are regularly reminded that preliminary assignments are intended to be straightforward and short, and that the problems are only checked for completeness at the beginning of discussion.  

\subsection{Enforcing attendance in discussions}
\label{subsec:enforcing-attendance}
Many math courses at our institution have discussion sections with a TA, but often they are optional and/or consist solely of review of content or homework problems.  In our classroom model, students need to be aware that discussions are a crucial part of the course, and not simply for review.  In addition to offering a participation grade, we communicated this in several ways:

\begin{itemize}
\item 
briefly mentioning in each lecture something that would be covered in the upcoming discussion, e.g., ``you will prove this theorem in discussion on Thursday'';

\item 
describing preliminary problems explicitly as preparation for discussions, and checking them at the beginning of discussion;

\item 
stating in the syllabus that we ``reserve the right to deduct one additional full letter grade if you miss too many classes, or if sufficient participation is not demonstrated during problem sessions'' and emphasizing this at the course onset;

\item 
ensuring both TA and instructor were present during discussion days (in truth, this was done so we could more readily visit groups and answer student questions, but as a byproduct it conveyed an added importance of the discussion days); and

\item 
handing out weekly homework assignments on the same sheet as the week's discussion problems (even though homework assignments were also posted online, this convention helps convey that discussions were not an add-on to the course, but rather an integral part of it).  

\end{itemize}

In the end, both times we taught the course, discussion days had the same average attendence as lecture days.  


%% file: sec-sample-discussions.tex

\section{Sample discussion activities}
\label{sec:samplediscussions}

\subsection{Modular arithmetic discussion worksheet}
\label{subsec:sample-discussion-modular}


\begin{figure}[p]
	\begin{center}
		\fbox{
			\includegraphics[height=\textheight-38pt]{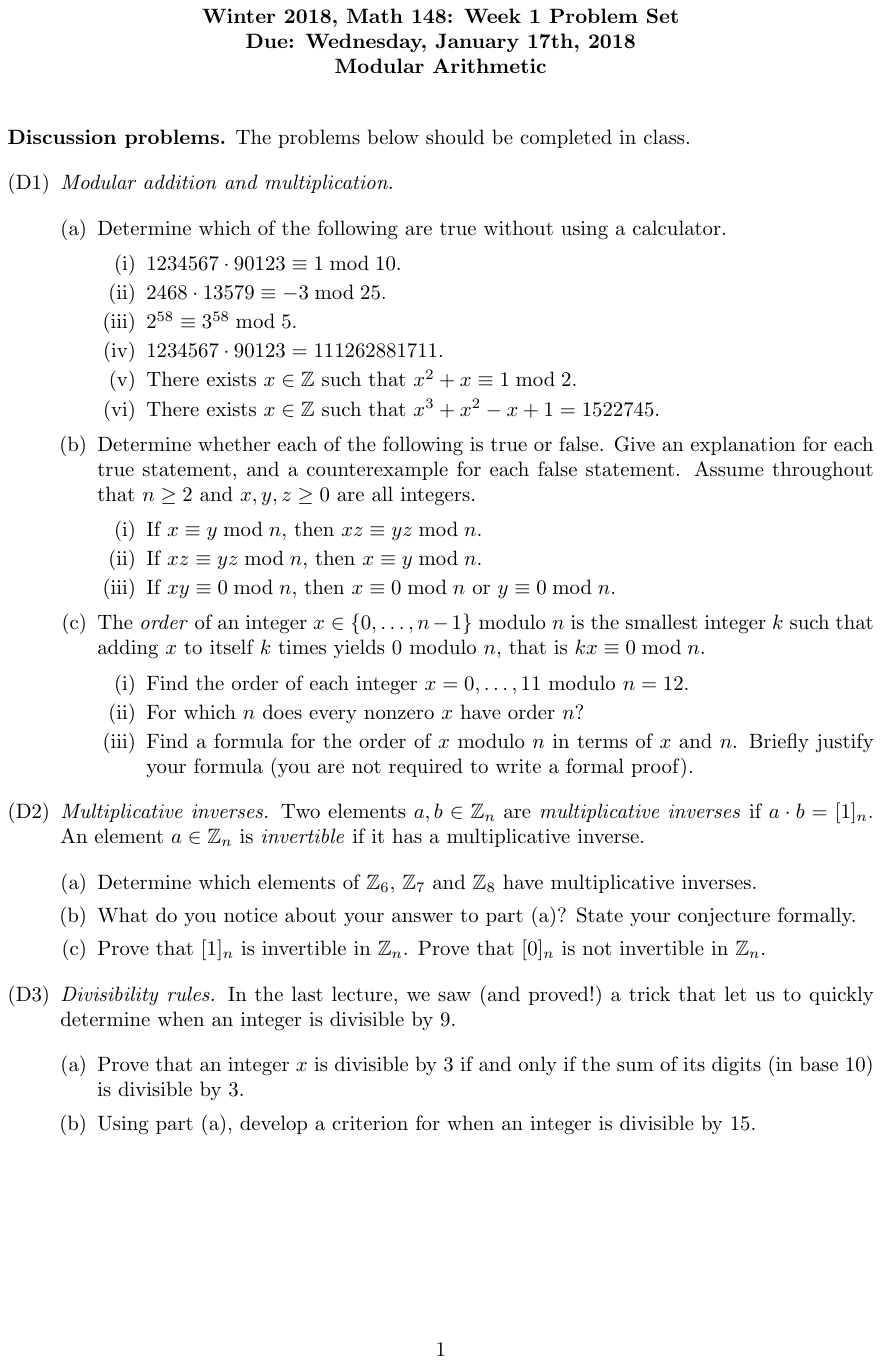}
		}
	\end{center}
	\caption{Week 1 discussion problems.}
	\label{f:d1disc}
\end{figure}

			

In Figure~\ref{f:d1disc}, we include the entire worksheet materials for the first week's activities.  
At this point in the course, students have only had one lecture on modular arithmetic, wherein emphasis was made on the fact that ``taking remainders'' commutes with addition and multiplication, and ``tricks'' for determining an integer's remainder modulo $9$ (summing the digits) and modulo $10$ (examining the last digit) were stated.  

Problem~(D1.a) is designed to help students employ modular arithmetic in some ways they have not yet seen, by gradually building on ideas they have.  Problems~(D1.a.i-iii) help them realize they can decrease their computation time by obtaining remainders before performing arithmetic, and develop a ``trick'' for finding remainders modulo~$25$.  For problems~(D1.a.iv) and~(D1.a.vi), students must go a step further and determine a clever way to use modular arithmetic (e.g., reducing to cases based on possible remainders), even though it is not explicitly present in the problem.  

Other discussion questions are intentionally open-ended, e.g.,\ problems~(D1.c.ii) and~(D1.c.iii),
which can inspire some creative answers from students, such as:
\begin{enumerate}[(D1.c.i)]
	\item[(D1.c.ii)] Answer: All $x$ such that the fraction $\nicefrac{x}{n}$ is reduced.
	
	\item[(D1.c.iii)] Answer: Reducing the fraction $\nicefrac{x}{n}$ to obtain $\nicefrac{p}{q}$, the order of $x$ mod $n$ is $q$.
	
\end{enumerate}
If we take these responses as evidence that the members of the group were more comfortable with reducing fractions than with terminology like \emph{relatively prime} and \emph{greatest common divisor}, then these answers were inherently more useful to that group. As instructors, this helped us more effectively define vocabulary and guide them through a proof that their answers were equivalent to those displayed on other groups' boards. 
It is also worth noting that \emph{order} is introduced here without appeal to groups. Later in the course, after presenting the definitions of \emph{group} and \emph{ring}, the general definition of order of a group element is introduced with familiar examples already at hand.

The Week 1 discussion sheet also introduces a major theme of this class, namely explicit computation of many small examples as a method for generating conjectures, e.g., problems~(D2.a) and~(D2.b).  
This process is repeated throughout the course, and we deliberately emphasize how many of the theorems ``handed down from on high'' in mathematics courses were in fact discovered through this kind of exploration. 

Although question (D2.b) does not ask for a proof of the supplied conjecture, this provides a natural place to extend the material for any groups proceeding through the worksheet quickly. Similarly, for the true/false statements from problem~(D1.b), a group could be asked to formally write a proof for the true statements, or asked if they could generalize a particular counterexample to a family of counterexamples and rigorously justify why every instance in the family would indeed be a counterexample.

\subsection{Block designs discussion activity}
\label{subsec:sample-discussion-block-designs}

\begin{figure}[t]
	\begin{center}
		\fbox{
			\includegraphics[width=\textwidth-0.4in]{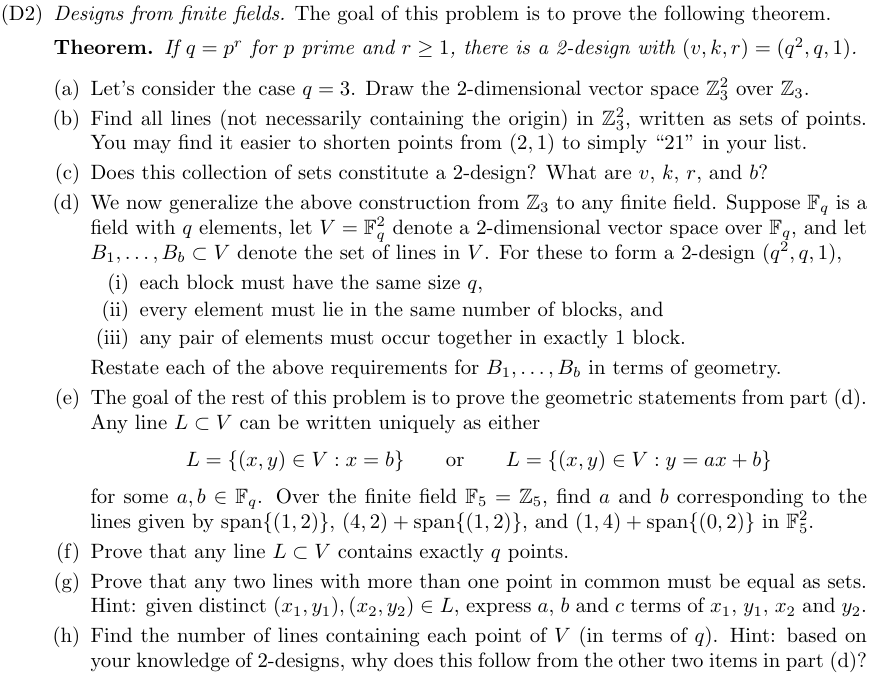}
		}
	\end{center}
	\caption{Week 7 discussion problem (D2).}
	\label{f:week7d2}
\end{figure}


Figure~\ref{f:week7d2} depicts a discussion problem that demonstrates the example-conjecture-proof appoach for guided discovery of ``big'' theorems or ideas.  In it, they are told up front the theorem they are to prove, and the ``discovery'' aspect lies in the main idea of the constructive proof:\ the blocks are lines in a finite vector space, and each requirement to be a 2-design can be interpreted as a geometric fact about such lines and their intersections.  

Students begin by working through the example $q = 3$.  The resulting 2-design has $12$ blocks, each with $3$ elements, and the corresponding vector space $\ZZ_3^2$ has $9$~points arranged in a $3 \times 3$ grid, making this example large enough to be instructive while not too large to be too time consuming or tedious.  The key connection to the geometry comes in part~(d), where students are asked to give a geometric interpretation of each requirement to be a 2-design.  The remainder of the activity guides students through proving each geometric statement obtained in part~(d), again motivated by specific examples that have been carefully chosen ahead of time to be an appropriate size and sufficient utility for identifying a general argument.  

This activity was to be completed in a $50$-minute class.  While not all groups completed the entire activity in that time, the most important take-aways are (i)~the formal statement of the main theorem, and (ii)~the geometric interpretations in part~(d), which students had no trouble completing within the allotted time.

\subsection{Projective geometry discussion activity}
\label{subsec:sample-discussion-projective}

\begin{figure}[t]
	\begin{center}
		\fbox{
			\includegraphics[width=\textwidth-0.4in]{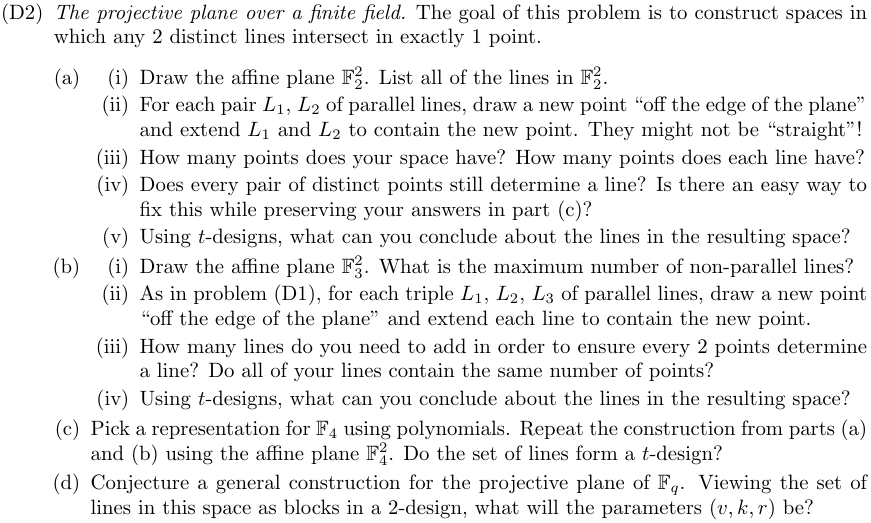}
		}
	\end{center}
	\caption{Week 8 discussion problem (D2).}
	\label{f:week8d2}
\end{figure}

\begin{figure}[t!]
	\begin{center}
		\begin{subfigure}[t]{0.28\textwidth}
			\begin{center}
				\includegraphics[width=\textwidth]{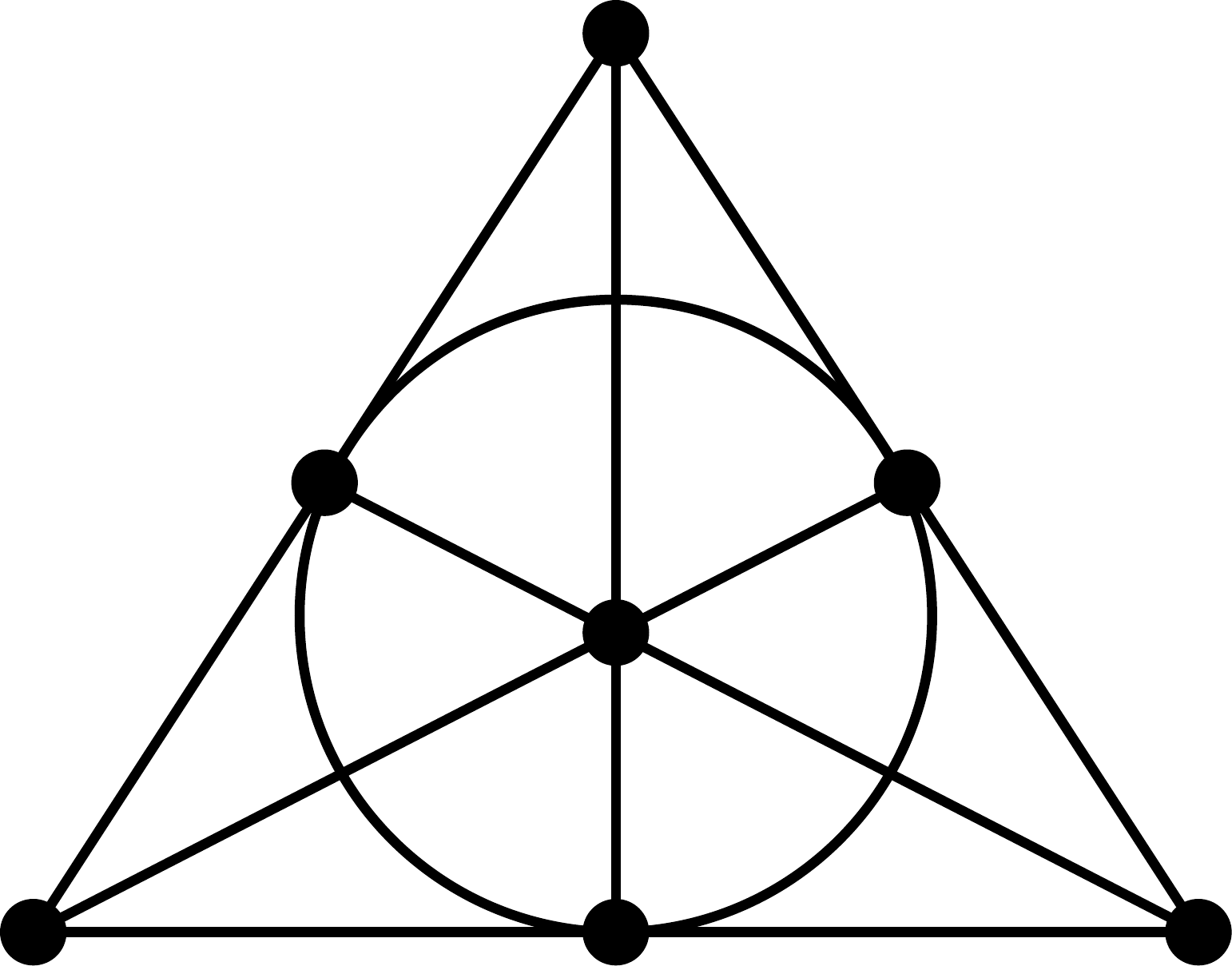}
			\end{center}
			\caption{}
			\label{f:rp2circle}
		\end{subfigure}
		\hspace{0.04\textwidth}
		\begin{subfigure}[t]{0.28\textwidth}
			\begin{center}
				\includegraphics[width=\textwidth]{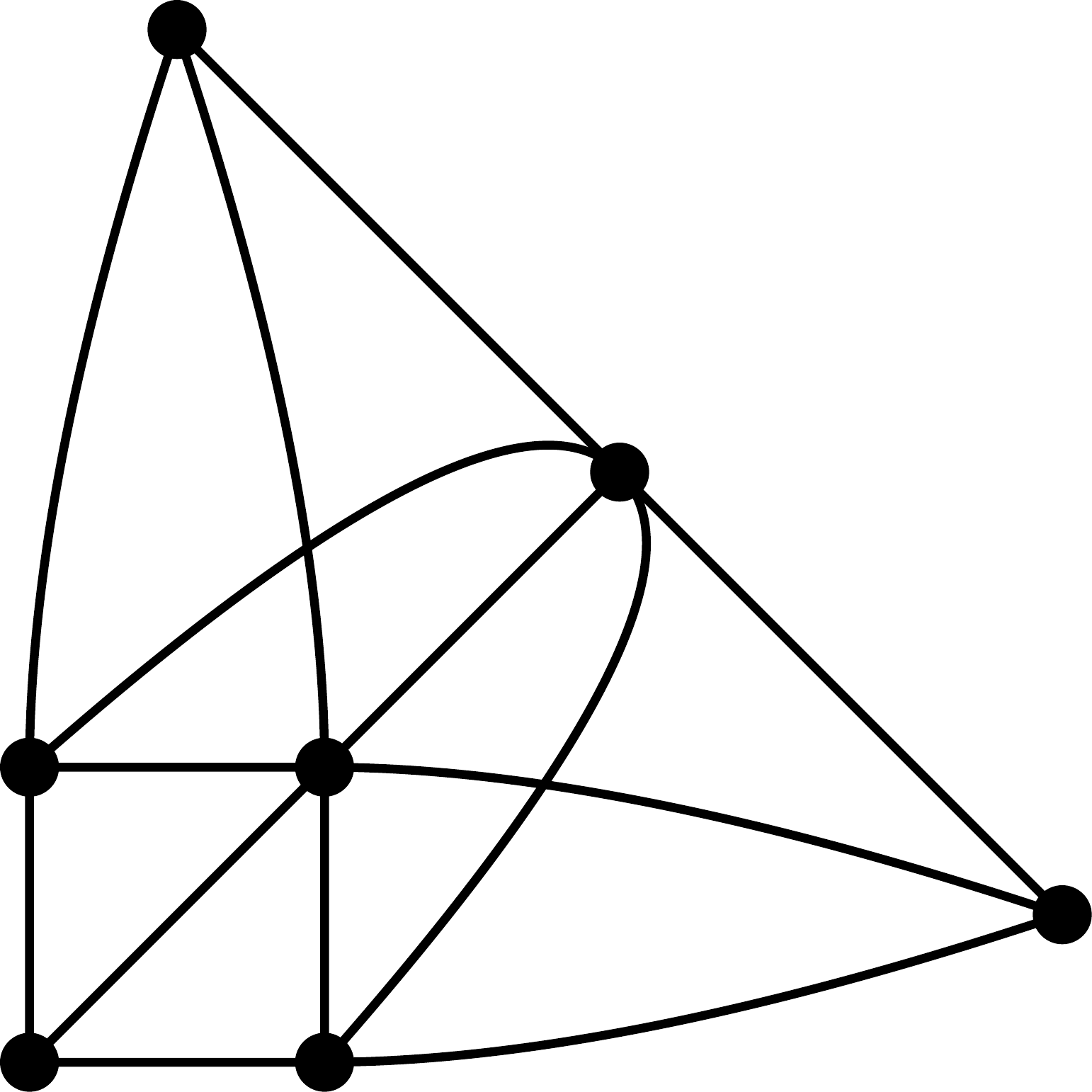}
			\end{center}
			\caption{}
			\label{f:rp2lineatinfinity}
		\end{subfigure}
		\hspace{0.04\textwidth}
		\begin{subfigure}[t]{0.28\textwidth}
			\begin{center}
				\includegraphics[width=\textwidth]{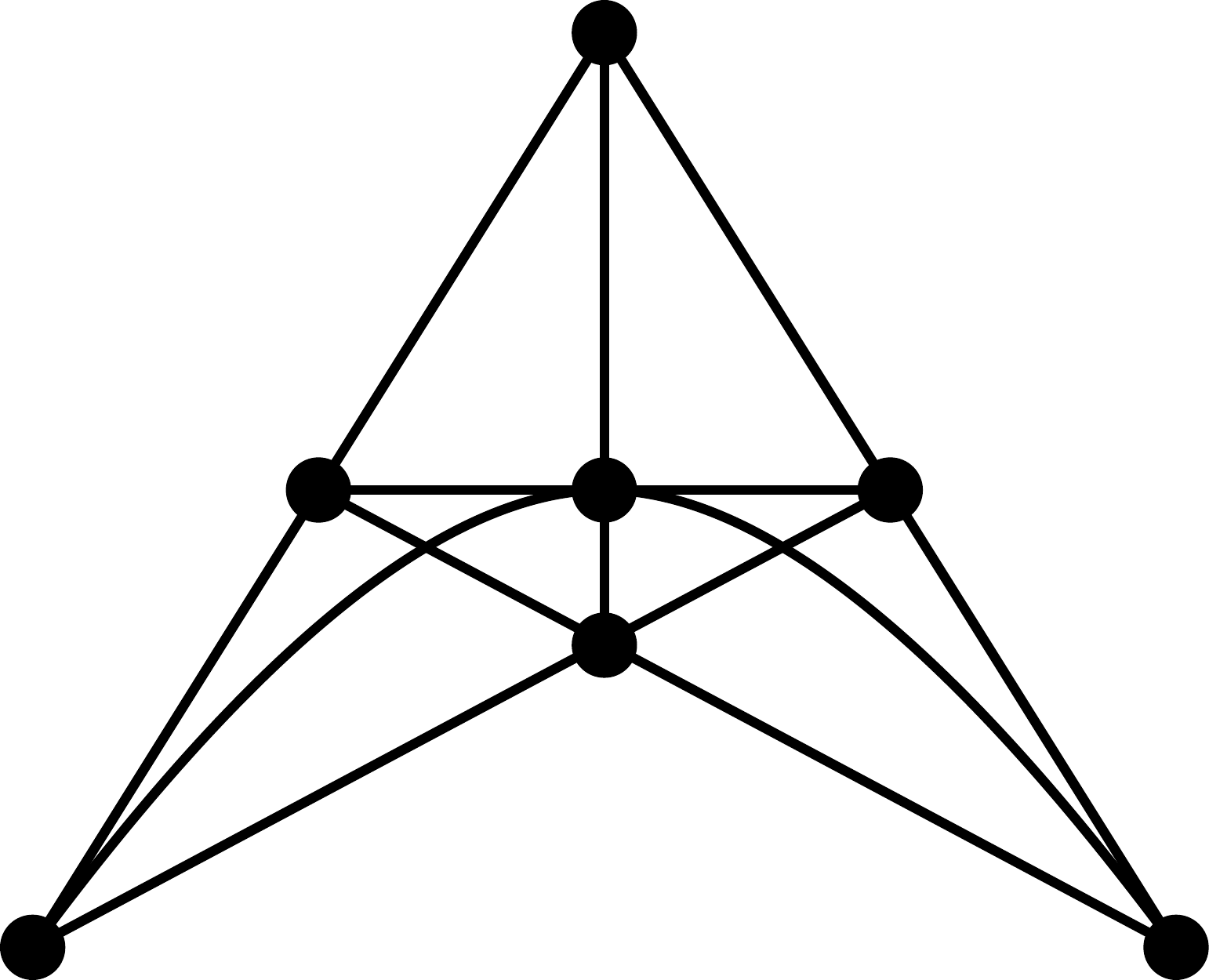}
			\end{center}
			\caption{}
			\label{f:rp2other}
		\end{subfigure}
	\end{center}
	\caption{Depictions of the projective plane over $\FF_2$ drawn in discussions.}
	\label{f:rp2}
\end{figure}

Discussion problem (D2) in Week 8, included here as Figure~\ref{f:week8d2}, introduces projective planes over finite fields.  Mathematically, the idea is to mimic the construction of the projective plane as the unioning of the affine plane $\RR^2$ and the line at infinity, but using $\FF_q$ instead of $\RR$.  The result is a finite set of points with a prescribed set of ``lines'' satisfying (i)~every line contains the same number of points, and (ii)~every pair of lines intersects exactly once.  

On the day of discussion, students are guided step by step to construct the projective plane over $\FF_2$, but not told initially what they are constructing (and most have not seen projective geometry before).  When students near the end of part~(a), we give the additional instruction to rearrange their drawing so the lines are ``as straight as possible''.  The goal of this intentionally ambiguous instruction is to give students an opportunity for open-ended exploration and to further internalize the geometric aspects of this construction.  As this discussion occurs late in the course, students are visibly more willing to ``just try things and see what happens'' without a concern for finding ``the" right answer, something several students commented on in their exit interviews.  

Surprisingly, most groups hit several of the same milestones.  Usually, the first historical drawing discovered (and arguably the most symmetric) is the one depicted in Figure~\ref{f:rp2circle}, wherein all but one line is straight.  When asked to further manipulate their drawing so that ``the least straight line is as straight as possible'', most groups discover that the further away they draw the line added in part (iv), the straighter the remaining lines become (Figure~\ref{f:rp2lineatinfinity}).  It is at this point that we point out this line is called the \emph{line at infinity} and represents where parallel lines meet at the horizon.  
Some groups do deviate from this, for instance settling on the arrangement depicted in Figure~\ref{f:rp2other} (which arguably best satisfies the ``lines as straight as possible'' among the three in Figure~\ref{f:rp2}).  This arrangment can be interpreted geometrically as having the line at infinity drawn horizontally through the middle.  

Around half way through class, students are instructed to move on to the remaining parts of problem~(D2).  Most important are part~(b), which repeats the construction over $\FF_3$, and part~(d), which ties these constructions back to the current topic of block designs.  Part~(c) is included to give extra practice with finite fields constructed as quotient rings, but was often skipped in the interest of time in order to get to part~(d).  

In~addition to the sense of discovery from the activity described above, this topic is a fantastic opportunity to introduce students to the broader subject of projective geometry, which, unlike some of the earlier material, is new to nearly every student.  The~following day's lecture includes a brief example-driven introduction to the projective plane over $\RR$ and its use in algebraic geometry, where the solution sets of some polynomial systems are more homogenous when viewed in projective space.  

When asked about the course, over half of the students listed this discussion day as their favorite, citing the geometric and pictorial aspects.  
The first time this course was taught, this discussion was the first one in which students worked at the board instead of at desks (see Section~\ref{subsec:discussion-changes}), which likely helped it make an extra impression.

%% file: sec-2ndtime.tex

\section{Changes the second time}
\label{sec:2ndtime}

\subsection{Content changes}
\label{subsec:content-changes}

Several topics and problem styles were more challenging for students than we originally anticipated, in particular those that did not require much background knowledge (hence, easy from the instructors' point of view), but did require more abstract thinking and reasoning (hence, hard from the students' point of view).

\subsubsection{Combinatorial reasoning}
In the first year, several discussion/homework problems in Week~5 (Applications of Finite Fields) had a combinatorial flavor.  

\begin{itemize}
\item
For $p$ prime, how many 1-dimensional linear subspaces does $\FF_p^2$ have?  

\item
How many 2-dimensional linear subspaces does $\FF_5^3$ have?  

\item
Find the number of bases of $\FF_{p^r}^2$, a 2-dimensional vector space over $\FF_{p^r}$.  

\end{itemize}

For the first question, for example, we expected students to come up with an argument similar to the following:

\begin{quote}
{A 1-dimensional linear subspace is a line through the origin. Every nonzero point in $\FF_p^2$ uniquely determines such a line. And every line in $\FF_p^2$ contains exactly $p - 1$ nonzero points. So the number of nonzero points, $p^2 - 1$, divided by the number of nonzero points per line, $p - 1$, gives the number of unique 1-dimensional subspaces, $(p^2 - 1)/(p - 1) = p + 1$.}
\end{quote}

At first, we thought of this argument as containing only one non-trivial statement:\ the span of a nonzero vector in $\FF_p^2$ always contains $p-1$ nonzero points. Since a previous question on the worksheet was to show this, we expected the problem to be straightforward. 
On closer examination, however, this proof requires some combinatorial sophistication. Counting lines by counting points that uniquely define them, and counting points by over-counting and then dividing, are both nontrivial leaps for students without practice at this sort of reasoning. We found it difficult to guide students to discover these strategies on their own. In~the second iteration of the course, we omitted these particular questions.

\subsubsection{Linear algebra review}
The second half of the course uses basic concepts from linear algebra to develop some families of error-correcting codes, which is why one quarter of linear algebra is a prerequisite for the course. In our first winter, we expected no difficulties using notions like basis, dimension, and subspace in our discussion of vector spaces over finite fields. However, a first course in linear algebra typically focuses on matrix operations, like row reduction and computing determinants. To these students, dimension might mean ``the number of pivot variables" more than it means ``the maximum number of linearly independent elements of a vector space."  Moreover, knowing the definition of dimension for real and complex vector spaces doesn't mean a student will automatically abstract this definition to more general fields (in this case, finite fields).  Although we only use linear algebra concepts that were covered in the prerequisite course, in the next course we spent more time reviewing linear algebra concepts and emphasizing how they translate to vector spaces over finite fields.

\subsubsection{Omitting some longer proofs from the lectures.}
The first time we taught the course, one or two lectures were spent proving the fundamental theorem of finite fields.  Our goal was to illustrate how to work with finite fields in a mathematically rigorous way, but this was a challenge, since (i)~most proofs require the use of technical machinery beyond the scope of this course (e.g., Galois groups), and (ii)~the disproportionate length and level of rigor felt out of place with the rest of the course.  The second time we taught the course, we opted for additional examples and discussion days intended to convey the intuition behind the fundamental theorem, rather than covering its mathematically rigorous development.

\subsection{Discussion changes}
\label{subsec:discussion-changes}

\subsubsection{All group work done on the board}
\
During our first time teaching this course, discussion groups mainly worked at their desks, which were rearranged into small clusters at the start of class. It wasn't until Week 8, when we introduced the projective plane over finite fields (see Section \ref{subsec:sample-discussion-projective}), that we took advantage of the many chalkboards in our classroom by having each group work at one.  In addition to being more fun, we quickly discovered that this change enhances discussions in multiple ways:
\begin{itemize}

\item 
It is easier to see what everyone is writing, compared to everyone working on paper at a cluster of desks, and the visibility creates some accountability. Groups are more attentive to what each member is writing, and spend more time coming to an agreement about what to display on the board.

\item 
Groups are forced to cooperate more. They may have to take turns with the chalk or markers, and edit and discuss each other's work. At the board, it's much harder for one student to move to a new question on their own, as they might do when sitting at a desk.

\item 
Taking on a physical role that is usually only performed by instructors reinforces student ownership of the material and the investigative nature of the discussion problems. A change in physical setting can also break students out of their usual group work patterns. 

\end{itemize}
It is worth pointing out that by having each group working on their own board space simultaneously, the social anxiety that might accompany working in front of a class is greatly reduced, while retaining the empowering and fun aspects of being at the board. 

This seemingly small change had such a positive response that we made it an integral aspect of the second course's discussions. We encourage anyone implementing a similar class to work with their college or university to schedule similar courses in rooms with adequate board space.

\subsubsection{Emphasizing that discussion worksheets did not need to be finished}

Even though discussion problems were never submitted from credit, some of them contained material that was relevant to homework and exam problems. 
Because the students didn't know \emph{a priori} which problems were relevant, they felt pressured to get through the entire discussion worksheet. 

The second time we taught the course, we addressed this issue by pointing discussion groups to the most essential parts of each worksheet. For example, we might make an announcement to the class after 25 of 50 minutes had passed:\ ``If any groups are still working on Discussion Problems 1 and 2, it's okay to leave those unfinished, but please move on to Problem 3. We want you to get through parts (a) and (b) of Problem 3, since these will be useful to you on the next homework assignment.''  By making such announcements partway through class, we helped avoid students skipping problems unnecessarily to do only the ones that ``really matter''.

%% file: sec-conclusion.tex
\section{Conclusion}
\label{sec:conclusion}

The authors did not conduct a rigorous study of student outcomes under this course design. While we were pleased by student evaluations and feedback---including anonymous comments such as ``I loved the structure of the class", "I actually wanted to go to class", and "My favorite math course at [the university]"---we have no direct comparisons of evaluation scores or student performances to lecture-only instances of the course. Our intention with this article is to offer our own reflections on the experience.


As instructors, we found that spending half of instruction time on small group work made the classroom more fun and engaging. We got to know our students better. We enjoyed having time to give interactive guidance and feedback in the classroom instead of relegating this to office hours only. When students had trouble with a definition, example, theorem, etc., we were able to identify the issue right away in class, instead of sending students home confused. And we frequently experienced one of the greatest joys of teaching:\ witnessing the ``aha moment" when a student suddenly figures out a mathematical concept on her own. 

%
%
%
%
%

Based on our experiences, both authors have continued to use this classroom structure at their respective universities, for a variety of classes including introductory number theory, introductory combinatorics, and undergraduate and graduate abstract algebra. Again, we have no rigorous analysis of specific student outcomes, but we informally observe that this approach is highly effective as well as fun. For other educators who find the evidence supporting active learning compelling, but find the thought of a flipped classroom daunting, we suggest giving the ``half lecture" method a try.